\documentclass[a4paper,12pt]{amsart}

\usepackage{amsmath,graphics}
\usepackage{amssymb}
\usepackage{amsfonts}
\usepackage{latexsym}
\usepackage{eucal}
\usepackage[dvips]{graphicx}
\usepackage{color}
\usepackage{enumerate}
\usepackage{mathrsfs}
\newtheorem{thm}{Theorem}[section]

\newtheorem{propo}[thm]{Proposition}
\newtheorem{lem}[thm]{Lemma}

\renewcommand{\Re}{{\rm Re}}
\renewcommand{\Im}{{\rm Im}}
\newcommand{\R}{\mathbb{R}}
\newcommand{\C}{\mathbb{C}}
\newcommand{\Z}{\mathbb{Z}}
\newcommand{\N}{\mathbb{N}}

\newcommand{\D}{ \mathcal{D}}

\newcommand{\hh}{\mathbb{H}^{2}}

\newcommand{\G}{{\bf G}}

\newcommand{\lt}{{\mathcal L}}

\newcommand{\Fp}{{\mathbb F}_p}
\newcommand{\half}{{\textstyle{\frac{1}{2}}}}
\newcommand{\rr}{{ \varrho}}

\begin{document}
\bibliographystyle{plain}
\title[L-functions and infinite index subgroups]{L-functions and sharp resonances of infinite index congruence subgroups of $SL_2(\Z)$}

\author[Dmitry Jakobson]{Dmitry Jakobson}
\address{McGill University\\
Department of Mathematics and Statistics\\
805 Sherbrooke Street West\\
Montreal, Quebec, Canada H3A0B9 
}
\email{jakobson@math.mcgill.ca}
\author{Fr\'ed\'eric Naud}
\address{Laboratoire de Math\'ematiques d'Avignon \\
Campus Jean-Henri Fabre, 301 rue Baruch de Spinoza\\
84916 Avignon Cedex 9. }
\email{frederic.naud@univ-avignon.fr}
 \maketitle
\begin{abstract} Let $\Gamma$ be a convex co-compact subgroup of $SL_2(\Z)$ and consider the "congruence subgroups"
$\Gamma(p)\vartriangleleft \Gamma$, for $p$ prime. Let $X_p:=\Gamma(p)\backslash \hh$ be the associated family of hyperbolic surfaces
covering $X:=\Gamma\backslash \hh$, we investigate the behaviour of resonances of the Laplacian on $X_p$ as $p$ goes to infinity. We prove a factorization formula for the Selberg zeta function $Z_{\Gamma(p)}(s)$ in term of $L$-functions $L_\Gamma(s,\rr)$ related to irreducible
representations $\rr$ of the Galois group $\G=SL_2(\Fp)$ of the covering, together with a priori bounds and analytic continuation. We use this factorization
property combined with an averaging technique over representations to prove a new existence result of non-trivial resonances in an effective low frequency strip. 
 \end{abstract}
 \tableofcontents
 \section{Introduction and results}
 In mathematical physics, resonances generalize the $L^2$-eigenvalues in situations where the underlying geometry is non-compact. Indeed, when the geometry has infinite volume, the $L^2$-spectrum of the Laplacian is mostly continuous and the natural replacement data for the missing eigenvalues are provided by resonances which arise from a meromorphic continuation of the resolvent of the Laplacian. 
 
To be more specific, in this paper we will work with the positive Laplacian $\Delta_X$ on hyperbolic surfaces $X=\Gamma \backslash \hh$, where $\Gamma$
 is a convex co-compact subgroup of $PSL_2(\R)$. A good reference on the subject is the book of Borthwick \cite{Borthwick}. Here $\hh$ is the hyperbolic plane endowed with its metric of constant curvature $-1$.
Let $\Gamma$ be a geometrically finite Fuchsian group of isometries acting on $\hh$. This means
that $\Gamma$ admits a finite sided polygonal fundamental domain in $\hh$. We will require that $\Gamma$ has no {\it elliptic} elements different from the identity and that the quotient $\Gamma \backslash \hh$ is of {\it infinite hyperbolic area}. We assume in addition in this
paper that $\Gamma$ has no parabolic elements (no cusps).
Under these assumptions, the quotient space 
$X=\Gamma \backslash \hh$ is a Riemann surface (called convex co-compact) whose {\it ends geometry} is as follows.
The surface $X$ can be decomposed into a compact surface $N$ with geodesic boundary, called the Nielsen region, on which  infinite area ends $F_i$ are glued : the funnels. 
A funnel $F_i$ is a half cylinder isometric to 
$$F_i=(\R /l_i \Z)_\theta \times (\R^+)_t,$$ where $l_i>0$, with the warped metric 
$$ds^2=dt^2+\cosh^2(t)d\theta^2.$$
The limit set $\Lambda(\Gamma)$ is defined as 
$$\Lambda(\Gamma):=\overline{\Gamma.z}\cap \partial \hh,$$
where $z\in \hh$ is a given point and $\Gamma.z$ is the orbit under the action of $\Gamma$ which accumulates
on the boundary $\partial \hh$. The limit set $\Lambda$ does not depend on the choice of $z$ and its Hausdorff dimension $\delta(\Gamma)$
is the critical exponent of Poincar\'e series \cite{Patterson}. 

\bigskip
The spectrum of $\Delta_X$ on $L^2(X)$ has been described completely by 
Lax and Phillips and Patterson in \cite{LP1,Patterson} as follows: 
\begin{itemize}
\item The half line $[1/4, +\infty)$ is the continuous spectrum.
\item There are no no embedded eigenvalues inside $[1/4,+\infty)$.
\item The pure point spectrum is empty if $\delta\leq \half$, and finite and starting at $\delta(1-\delta)$ if $\delta>\half$.
\end{itemize}

Using the above notations, the resolvent 
$$R(s):=(\Delta_X-s(1-s) )^{-1}:L^2(X)\rightarrow L^2(X)$$
is a holomorphic family for $\Re(s) >\half$, except at a finite number of possible poles related to the eigenvalues. From the work of Mazzeo-Melrose \cite{MM}, it can be meromorphically continued (to all $\C$)  from $C_0^\infty(X)\rightarrow C^\infty(X)$, and poles are called {\it resonances}. We denote
in the sequel by $\mathcal{R}_X$ the set of resonances, written with multiplicities.

\bigskip
To each resonance $s\in \C$ (depending on multiplicity) are associated generalized eigenfunctions (so-called purely outgoing states) $\psi_s \in C^\infty(X)$
which provide stationary solutions of the automorphic {\it wave equation} given by 
$$\phi(t,x)=e^{(s-\half)t}\psi_s(x),$$
$$\left (D_t^2+\Delta_X-\frac{1}{4} \right)\phi=0.$$
From a physical point of view, $\Re(s)-\half$ is therefore a rate of decay while $\Im(s)$ is a frequency of oscillation. Resonances that live the longest are called {\it sharp resonances} and are those for which $\Re(s)$ is the closest to the unitary axis $\Re(s)=\half$. In general, $s=\delta$ is the only explicitly known resonance (or eigenvalue if $\delta>\half$). There are very few effective results on the existence of {\it non-trivial } sharp resonances,
and to our knowledge the best statement so far is due to the authors \cite{JN2}, where it is proved that for all $\epsilon>0$, there are infinitely many resonances in the strip
$$\left \{ \Re(s) >\frac{\delta(1-2\delta)}{2}-\epsilon  \right \}.$$
It is conjectured in the same paper \cite{JN2} that for all $\epsilon>0$, there are infinitely many resonances in the strip $\{ \Re(s)>\delta/2-\epsilon\}$. However, the above result, while proving existence of non-trivial resonances, does not provide estimates on the imaginary parts (the frequencies), and it is a notoriously hard problem to locate precisely non trivial resonances.
The goal of the present work is to obtain a different type of existence result by looking at families of "congruence" surfaces.  
Let $\Gamma$ be an infinite index, finitely generated, free subgroup of $SL_2(\Z)$, without parabolic elements. Because $\Gamma$ is free, the projection map
$\mathcal{\pi}:SL_2(\R)\rightarrow PSL_2(\R)$ is injective when restricted to $\Gamma$ and we will thus identify $\Gamma$ with $\pi(\Gamma)$, i.e.
with its realization as a Fuchsian group. Under the above hypotheses, $\Gamma$ is a convex co-compact group of isometries. For all $p>2$ a prime number, we define the congruence subgroup $\Gamma(p)$ by 
$$\Gamma(p):=\{ \gamma \in \Gamma \ :\ \gamma \equiv \mathrm{Id}\ \mathrm{mod}\ p\},$$
and we set $\Gamma(0)=\Gamma$.
Recently, these "infinite index congruence subgroups" have attracted a lot of attention because of the key role they play in number theory and graph theory. We mention the early work of Gamburd \cite{Gamburd1} and the more recent by Bourgain-Gamburd-Sarnak \cite{BGS}, Bourgain-Kontorovich \cite{Kontorovich1} and Oh-Winter 
\cite{OhWinter}. In all of the previously mentioned works, the spectral theory of surfaces
$$X_p:=\Gamma(p)\backslash \hh,$$
plays a critical role and knowledge on resonances is mandatory, in particular knowledge on uniform spectral gaps as $p\rightarrow +\infty$.
In \cite{JN1}, the authors have started investigating the behaviour of resonances 
in {\it the large $p$ limit} and the present paper goes in the same direction with different techniques involving sharper tools of representation theory.

A way to attack any problem on resonances of hyperbolic surfaces is through the Selberg zeta function defined for $\Re(s)>\delta$ by
$$Z_\Gamma(s):=\prod_{\mathcal{C}\in \mathcal{P}} \prod_{k\in \N} \left (1-e^{-(s+k)l(\mathcal{C})} \right), $$
where $\mathcal{P}$ is the set of primitive closed geodesics and $l(\mathcal{C})$ is the length. This zeta function extends analytically to $\C$ and it is known from the work of Patterson-Perry \cite{PatPerry} that non-trivial zeros of $Z_\Gamma(s)$ are resonances with multiplicities. Our first goal is to establish a factorization formula for $Z_{\Gamma (p)}(s)$ using the representation theory of $SL_2(\Fp)$. Indeed, it is known from Gamburd \cite{Gamburd1}, that the map 
$$\pi_p:\left \{ \Gamma \rightarrow  SL_2(\Fp) \atop \gamma \mapsto \gamma\ \mathrm{mod}\ p \right.$$
is onto for all $p$ large, and we thus have a family of Galois covers $X_p\rightarrow X$ with Galois group $SL_2(\Fp)$. Let $\{ \rr \}$ denote the set
of irreducible complex representations of $\G:=SL_2(\Fp)$, and given $\rr$ we denote by $\chi_{\rr}$ its character, $V_\rr$ its linear representation space and we set 
$$d_\rr:=\mathrm{dim}(V_\rr).$$ Our first result is the following.
\begin{thm}
\label{main1} 
For $\Re(s)>\delta$, consider the L-function defined by
$$L_\Gamma(s,\rr):=\prod_{\mathcal{C}\in \mathcal{P}} \prod_{k\in \N} \det \left (Id_{V_\rr}-\rr(\mathcal{C^k})e^{-(s+k)l(\mathcal{C})} \right),$$
where $\rr(\mathcal{C})$ is understood as $\rr(\pi_p(\gamma_{\mathcal{C}}))$ where $\gamma_{\mathcal{C}}\in \Gamma$ is any representative
of the conjugacy class defined by $\mathcal{C}$. Then we have the following facts.
\begin{enumerate}
 \item For all $\rr$ irreducible, $L_\Gamma(s,\rr)$ extends as analytic function to $\C$.
 \item There exist $C_1,C_2>0$ such that for all $p$ large, all $\rr$ irreducible representation of $\G$, and all $s\in \C$,
 we have
 $$\vert L_\Gamma(s,\rr) \vert \leq C_1 \exp\left(C_2 d_\rr \log(1+d_\rr)(1+ \vert s\vert^2) \right).$$
 \item For all $p$ large, we have the formula
 $$Z_{\Gamma(p)}(s)=\prod_{\rr\  \mathrm{irreducible}} \left(L_\Gamma(s,\rr) \right)^{d_\rr}.$$
\end{enumerate}
\end{thm}
Notice that the $L$-function for the {\it trivial representation} is just $Z_\Gamma(s)$ and thus $Z_\Gamma(s)$ is always a factor of
$Z_{\Gamma(p)}(s)$.  There is a long story of L-functions associated with compact extensions of geodesic flows in negative curvature, see for example \cite{Sarnak,KatSun} and \cite{ParryPollicott1}. In the case of pairs of Hyperbolic
pants with symmetries, a similar type of factorization has been considered for numerical purposes by Borthwick and Weich \cite{BW1}.
The above factorization is very similar to the factorization of Dedekind zeta functions as a product of Artin L-functions in the case of Number fields. In the context of hyperbolic surfaces with infinite volume, the above statement is new and interesting in itself for various applications. 

\bigskip
In \cite{JN1}, by combining trace formulae techniques with some a priori upper bounds for $Z_{\Gamma(p)}(s)$ obtained via transfer operator techniques,
we proved the following fact. For all $\epsilon>0$, there exists $C_\epsilon>0$ such that for all $p$ large enough,  
$$C_\epsilon^{-1} p^3\leq \#\mathcal{R}_{X_p}\cap \{ \vert s\vert \leq (\log(p))^\epsilon\}\leq C_\epsilon p^3 (\log(p))^{1+2\epsilon}.$$
We point out that $p^3\asymp \mathrm{Vol}(N_p)$, where $\mathrm{Vol}(N_p)$ is the volume of the convex core of $X_p$, therefore these bounds can be thought as a {\it Weyl law} in the large $p$ regime. In the case of covers of compact or finite volume manifolds, precise results for the Laplace spectrum in the "large degree" limit have been obtained in the past
in \cite{Degeorge, Donnelly}. In the case of infinite volume hyperbolic manifolds, we also mention the density bound obtained by Oh \cite{Oh1}.

While this result has near optimal upper and lower bounds, it does not provide a lot of information
on the precise location of this wealth of non trivial resonances. The main result of this paper is as follows.

\begin{thm}
\label{main2}
Using the above notations, assume that $\delta>\frac{3}{4}$. Then for all $\epsilon,\beta>0$, and
for all $p$ large,
$$\# \mathcal{R}_{X_p}\cap \left \{ \delta-\textstyle{\frac{3}{4}}-\epsilon\leq \Re(s)\leq \delta\ \mathrm{and}\ 
\vert \Im(s) \vert \leq \left (\log(\log(p))\right)^{1+\beta}\right \} \geq \frac{p-1}{2}.$$
\end{thm}
Existence of convex co-compact subgroups $\Gamma$ of $SL_2(\Z)$ with $\delta_\Gamma$ arbitrarily close to $1$ is guaranteed by a theorem
of Lewis Bowen \cite{LewisBowen}. See also \cite{Gamburd1} for some hand-made examples. 
Theorem \ref{main2} shows that as $p\rightarrow \infty$, there are plenty of resonances with $\Re(s)\geq \delta-\textstyle{\frac{3}{4}}-\epsilon$ and small
imaginary part $\Im(s)$, providing an {\it existence theorem} with an effective control of both decay rate and frequency. Notice that from \cite{JN2}, we also know  that
for any $\sigma\leq \delta$, $\beta>0$, and all large $p$, $$\# \mathcal{R}_{X_p}\cap \left \{ \sigma \leq \Re(s)\leq \delta\ \mathrm{and}\ 
\vert \Im(s) \vert \leq \left (\log(\log(p))\right)^{1+\beta}\right \}$$
$$\leq C_{\beta,\sigma} \left (p^3 \log(p)\right).$$
One of the points of Theorem \ref{main2} is that we can produce non-trivial resonances without having to affect $\delta$, just by moving to a finite cover. The outline of the proof is as follows. Having established the factorization formula, we first notice that since the dimension of any non trivial representation of $\G$ is at least $\textstyle{\frac{p-1}{2}}$, it is enough to show that at least one of the $L$-functions
$L_\Gamma(s,\rr)$ vanishes in the described region as $p\rightarrow \infty$. We achieve this goal through an averaging technique (over irreducible $\rr$)
which takes in account the "explicit" knowledge of the conjugacy classes of $\G$, 
together with the high multiplicities in the length spectrum of $X$. Unlike in finite volume cases where one can 
take advantage of a precise location of the spectrum (for example by assuming GRH), none of this strategy applies here which makes it much harder to mimic existing techniques from analytic number theory.

The remainder of the paper is divided in $3$ sections. In $\S 2$ we prove analytic continuation of $L$-functions together with an a priori bound using transfer operators and singular value estimates.
In $\S 3$, we relate the existence of zero-free regions for $L$-functions to upper bounds on certain twisted sums over closed geodesics via some arguments of harmonic and complex analysis. In $\S 4$, we recall the algebraic facts around the group $\G=SL_2(\Fp)$ and apply it to prove Theorem \ref{main2} via averaging over $\rr$ to produce a relevant lower bound and conclude by contradiction.

\bigskip
{\noindent \bf Acknowledgements}. Both Authors are supported by ANR grant "GeRaSic". DJ was partially supported by NSERC, FRQNT and Peter Redpath fellowship. FN is supported by Institut Universitaire de France.

\section{Vector valued transfer operators and analytic continuation}

\subsection{Bowen-Series coding and transfer operator}
The goal of this section is to prove Theorem \ref{main1}. The technique follows closely previous works \cite{Naud1,JN1} with the notable addition that
we have to deal with vector valued transfer operators. We start by recalling Bowen-Series coding and holomorphic function spaces needed for our analysis. Let $\hh$ denote the Poincar\'e upper half-plane 
$$\hh=\{ x+iy\in \C\ :\ y>0\}$$
endowed with its standard metric of constant curvature $-1$
$$ds^2=\frac{dx^2+dy^2}{y^2}.$$ 
The group of isometries of $\hh$ is $\mathrm{PSL}_2(\R)$ through the action of 
$2\times 2$ matrices viewed as M\"obius transforms
$$z\mapsto \frac{az+b}{cz+d},\ ad-bc=1.$$ 
Below we recall the definition of Fuchsian Schottky groups which will be used to define transfer operators.
A Fuchsian Schottky group is a free subgroup of $\mathrm{PSL}_2(\R)$ built as follows. Let $\D_1,\ldots, \D_m,\D_{m+1},\ldots, \D_{2m}$, $m\geq 2$, 
be $2m$ Euclidean {\it open} discs in $\C$ orthogonal to the line $\R\simeq \partial \hh$. We assume that for all $i\neq j$, $\overline{\D_i} \cap \overline{\D_j}=\emptyset$. 
Let $\gamma_1,\ldots,\gamma_m \in \mathrm{PSL}_2(\R)$ be $m$ isometries such that for all $i=1,\ldots,m$, we have
$$\gamma_i(\D_i)=\widehat{\C}\setminus \overline{\D_{m+i}},$$
where $\widehat{\C}:=\C\cup \{ \infty \}$ stands for the Riemann sphere. For notational purposes, we also set $\gamma_i^{-1}=:\gamma_{m+i}$.
\begin{center}
\includegraphics[scale=0.65]{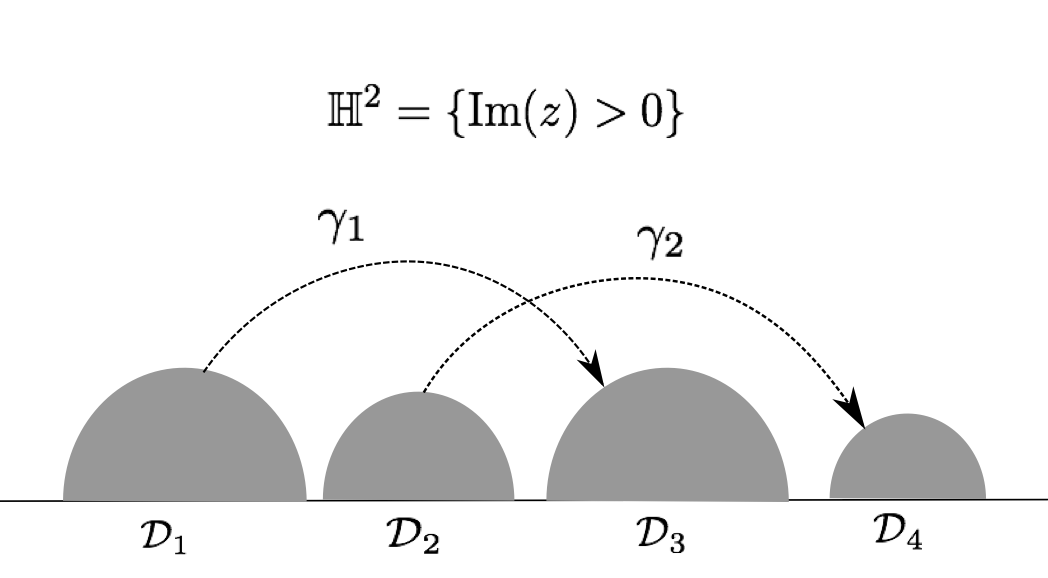}
\end{center}

\bigskip \noindent
Let $\Gamma$ be the free group generated by $\gamma_i,\gamma_i^{-1}$ for $i=1,\ldots,m$, then $\Gamma$ is a convex co-compact group, i.e. it is finitely generated and has no non-trivial parabolic element. The {\it converse is true} : up to isometry, convex co-compact hyperbolic surfaces
can be obtained as a quotient by a group as above, see \cite{Button}.

We assume from now on that each $\gamma_i \in \mathrm{PSL}_2(\R)$ comes from an element in $SL_2(\Z)$ so that $\Gamma$ is naturally identified with an infinite index, finitely generated free subgroup of $SL_2(\Z)$. 

For all $j=1,\ldots,2m$, set $I_j:=\D_j\cap \R$. One can define a map 
$$T:I:=\cup_{j=1}^{2m}I_j\rightarrow\R\cup\{\infty\}$$
by setting
$$T(x)=\gamma_j(x)\ \mathrm{if}\ x\in I_j.$$
This map encodes the dynamics of the full group $\Gamma$, and is called the Bowen-Series map, see 
\cite{BowenSeries1, Bowen1} for the genesis of these type of coding. The key properties are orbit equivalence and
uniform expansion of $T$ on the maximal invariant subset $\cap_{n\geq 1} T^{-n}(I)$ which coincides with the limit set
$\Lambda(\Gamma)$, see \cite{Bowen1}.

We now define the function space and the associated transfer operators.
Set 
$$\Omega:=\cup_{j=1}^{2m}\D_j.$$
 Each complex representation space $V_\rr$ is endowed with an inner product $\langle .,.\rangle_\rr$ which makes each representation
$$\rr: \G\rightarrow \mathrm{End}(V_\rr)$$
unitary. Consider now the Hilbert space $H^2_\rr(\Omega)$ which is defined as the set of {\it vector valued holomorphic functions} 
$F:\Omega\rightarrow V_\rr$ such that
$$\Vert F\Vert_{H^2_\rr}^2:=\int_{\Omega} \Vert F(z)\Vert^2_\rr dm(z)<+\infty, $$
where $dm$ is Lebesgue measure on $\C$. On the space $H^2_\rr(\Omega)$, we define a "twisted" by $\rr$ transfer operator $\lt_{\rr,s}$ by
$$\lt_{\rr,s}(F)(z):=\sum_{j} ((T')(T_j^{-1}))^{-s}F(y)\rr(T_j^{-1})=
\sum_{j\neq i} (\gamma_j')^s F(\gamma_j z)\rr(\gamma_j),\ \mathrm{if}\ z\in \D_i,$$
where $s\in \C$ is the spectral parameter. Here $\rr(\gamma_j)$ is understood as
$$\rr(\pi_p(\gamma_j)),\ \gamma_j\in SL_2(\Z),$$
where $\pi_p:SL_2(\Z)\rightarrow SL_2(\Fp)$ is the natural projection.
We also point out that the linear map $\rr(g)$ acts "on the right" on vectors $U \in V_\rr$ simply by fixing an orthonormal basis $\mathcal{B}=(e_1,\ldots,e_{d_\rr})$ of $V_\rr$ and setting
$$U\rr(g):=(U_1,\ldots,U_{d_\rr}) {\mathcal Mat}_{\mathcal{B}}(\rho(g)).$$
Notice that for all $j\neq i$, $\gamma_j:\D_i\rightarrow \D_{m+j}$ is a holomorphic contraction since $\overline{\gamma_j(\D_i)}\subset \D_{m+j}$.
Therefore, $\lt_{\rr,s}$ is a compact {\it trace class} operator and thus has a {\it Fredholm determinant}. We start by recalling a few facts.

We need to introduce some more notations. Considering a finite sequence $\alpha$ with
\[\alpha=(\alpha_1,\ldots,\alpha_n)\in \{1,\ldots, 2m\}^n,\]
we set 
\[ \gamma_\alpha:=\gamma_{\alpha_1}\circ \ldots \circ \gamma_{\alpha_n}. \]
We then denote by $\mathscr{W}_n$ the set of admissible sequences of length $n$ by
\[ \mathscr{W}_n:=\left \{ \alpha \in \{1,\ldots, 2m\}^n\ :\ 
\forall\ i=1,\ldots,n-1,\ \alpha_{i+1}\neq \alpha_i +m\ \mathrm{mod}\ 2m \right \}.\]
The set $\mathscr{W}_n$ is simply the set of reduced words of length $n$.
For all $j=1,\ldots, 2m$, we define $\mathscr{W}_n^j$ by
\[ \mathscr{W}_n^j:=\{ \alpha \in \mathscr{W}_n\ :\ \alpha_n\neq j \}. \] 
If $\alpha \in \mathscr{W}_n^j$, then $\gamma_\alpha$ maps $\overline{\D_j}$ into $\D_{\alpha_1+m}$. Using this set of notations, we have the formula for all
$z\in \D_j$, $j=1,\ldots,2m$, 
$$\lt_{\rr,s}^N(F)(z)=\sum_{\alpha \in \mathscr{W}_N^j} (\gamma_\alpha'(z))^s F(\gamma_\alpha z)\rr(\gamma_\alpha).$$
A key property of the contraction maps $\gamma_\alpha$ is that they are {\it eventually uniformly contracting}, see \cite{Borthwick}, prop 15.4 : there exist $C>0$ and $0<\rho_2<\rho_1<1$ such that for all $z\in \D_j$, for all $\alpha \in \mathscr{W}_n^j$ we have for all $n\geq 1$, 
\begin{equation}
\label{Ucont}
C^{-1}\rho_2^N\leq \sup_{z\in \D_j} \vert \gamma_\alpha'(z) \vert\leq C \rho_1^n
\end{equation}
In addition, they have the {\it bounded distortion property} (see \cite{Naud1} for proofs):
There exists $M_1>0$ such that for all $n,j$ and all $\alpha \in \mathscr{W}_n^j$, we have for all $z \in \D_j$,
\begin{equation}
\label{Bdist}
\left \vert \frac{\gamma''_\alpha(z)}{\gamma'_\alpha(z)}\right \vert \leq M_1.
\end{equation}
We will also need to use the {\it topological pressure} as a way to estimate certain weighted sums over words. 
We will rely on the following fact \cite{Naud1}. Fix $\sigma_0 \in \R$, then there exists $C(\sigma_0)$ such that for all $n$ and 
$\sigma\geq \sigma_0$, we have
\begin{equation}
\label{Pressure0}
\sum_{j=1}^{2m} \left ( \sum_{\alpha \in \mathscr{W}_n^j} \sup_{\D_j} \vert \gamma'_\alpha \vert^\sigma \right) 
\leq C_0e^{nP(\sigma_0)}.
\end{equation}
Here $\sigma\mapsto P(\sigma)$ is the {\bf topological pressure}, which is a strictly convex decreasing function which vanishes at $\sigma=\delta$, see \cite{Bowen1}.  In particular, whenever $\sigma>\delta$, we have $P(\sigma)<0$.
A definition of $P(\sigma)$ is by a variational formula:
$$P(\sigma)=\sup_{\mu} \left ( h_\mu(T)-\sigma \int_{\Lambda}\log \vert T' \vert d\mu \right), $$
where $\mu$ ranges over the set of $T$-invariant probability measures, and $h_\mu(T)$ is the measure theoretic entropy.
For general facts on topological pressure and theormodynamical formalism we refer to \cite{ParryPollicott2}.
Since we will only use it once for the spectral radius estimate below, we don't feel the need to elaborate more on various other definitions of the topological pressure and refer the reader to the above references.
\bigskip
\subsection{Norm estimates and determinant identity}
We start with an a priori norm estimate that will be used later on, see also \cite{JN1} where a similar bound (on a different function space) is proved in appendix.
\begin{propo}
\label{norm1}
Fix $\sigma=\Re(s) \in \R$, then there exists $C_\sigma>0$, independent of $p,\rr$ such that for all $s \in \C$  with $\Re(s)=\sigma$ and all 
$N$ we have
$$\Vert \lt_{\rr,s}^N\Vert_{H^2_\rr}\leq C_\sigma e^{C_\sigma \vert \Im(s)\vert} e^{2NP(\sigma)}.$$
\end{propo}
\noindent {\it Proof}. First we need to be more specific about the complex powers involved here. First we point out that given $z\in \D_i$ then for all $j\neq i$, $\gamma'_j(z)$ belongs to $\C\setminus (-\infty,0]$, simply because each $\gamma_j$ is in $PSL_2(\R)$. 
This make it possible to define $\gamma'_j(z)^s$ by 
$$\gamma'_j(z)^s:=e^{s\mathbb{L}(\gamma'_j(z))},$$
where $\mathbb{L}(z)$ is the complex logarithm defined  on $\C\setminus (-\infty,0]$ by the contour integral
$$\mathbb{L}(z):=\int_1^z\frac{d\zeta}{\zeta}.$$
By analytic continuation, the same identity holds for iterates.
In particular, because of bound (\ref{Ucont}) and also bound (\ref{Bdist}) one can easily show that there exists $C_1>0$ such that for all $N,j$ and all 
$\alpha \in \mathscr{W}_N^j$, we have
\begin{equation}
\label{Bound1}
\sup_{z\in \D_j} \vert \gamma'_\alpha(z)^s\vert \leq e^{C_1\vert \Im(s)\vert}\sup_{\D_j}\vert \gamma'_\alpha \vert^\sigma,
\end{equation}
where $\sigma=\Re(s)$.
We can now compute, given $F\in H^2_\rr(\Omega)$,
$$\Vert \lt_{\rr,s}^N(F)\Vert^2_{H^2_\rr}:=\sum_{j=1}^{2m}\sum_{\alpha,\beta \in \mathscr{W}_N^j}
\int_{\D_j} \gamma'_\alpha(z)^s \overline{\gamma'_\beta(z)^s} 
\langle F(\gamma_\alpha z) \rr(\gamma_\alpha),F(\gamma_\beta z)\rr(\gamma_\beta) \rangle_\rr dm(z).$$
By unitarity of $\rr$ and Schwarz inequality we obtain
$$\Vert \lt_{\rr,s}^N(F)\Vert^2_{H^2_\rr}\leq  e^{2C_1\vert \Im(s)\vert}\sum_j \sum_{\alpha,\beta}
\sup_{\D_j}\vert \gamma'_\alpha \vert^\sigma \sup_{\D_j}\vert \gamma'_\beta \vert^\sigma
\int_{D_j} \Vert F(\gamma_\alpha z) \Vert_\rr \Vert F(\gamma_\beta z) \Vert_\rr dm(z).$$ 
We now remark that $z\mapsto F(z)$ has components in $H^2(\Omega)$, the Bergman space of $L^2$ holomorphic functions
on $\Omega=\cup_j \D_j$, so we can use the scalar reproducing kernel $B_{\Omega}(z,w)$ to write (in a vector valued way)
$$F(\gamma_\alpha z )=\int_{\Omega}F(w) B_{\Omega}(\gamma_\alpha z,w)dm(w).$$
Therefore we get
$$\Vert F(\gamma_\alpha z )\Vert_\rr\leq \int_{\Omega}\Vert F(w)\Vert_\rr \vert B_{\Omega}(\gamma_\alpha z,w)\vert dm(w), $$
and by Schwarz inequality we obtain
$$\sup_{z\in \D_j}  \Vert F(\gamma_\alpha z )\Vert_\rr\leq \Vert F \Vert_{H^2_\rr} 
\left (\int_{\Omega} \vert B_{\Omega}(\gamma_\alpha z,w)\vert^2 dm(w) \right )^{1/2}.$$
Observe now that by uniform contraction of branches $\gamma_\alpha:\D_j\rightarrow \Omega$,  there exists a compact
subset $K\subset \Omega$ such that for all $N,j$ and $\alpha\in \mathscr{W}_{N}^j$,
$$\gamma_\alpha(\D_j)\subset K.$$
We can therefore bound 
$$\int_{\Omega} \vert B_{\Omega}(\gamma_\alpha z,w)\vert^2 dm(w)\leq C$$
uniformly in $z,\alpha$. We have now reached
$$ \Vert \lt_{\rr,s}^N(F)\Vert^2_{H^2_\rr}\leq  \Vert F \Vert_{H^2_\rr}^2 C_2e^{2C_1\vert \Im(s)\vert}\sum_j \sum_{\alpha,\beta}
\sup_{\D_j}\vert \gamma'_\alpha \vert^\sigma \sup_{\D_j}\vert \gamma'_\beta \vert^\sigma,$$
and the proof is now done using the topological pressure estimate (\ref{Pressure0}). $\square$

\bigskip
The main point of the above estimate is to obtain a bound which is independent of $d_\rr$. In particular the spectral radius 
$\rho_{sp}(\lt_{\rr,s})$ of $\lt_{\rr,s}:H^2_\rr(\Omega)\rightarrow H^2_\rr(\Omega)$ is bounded by
\begin{equation}
\label{radius1}
\rho_{sp}(\lt_{\rr,s})\leq e^{P(\Re(s))},
\end{equation}
which is uniform with respect to the representation $\rr$, and also shows that it is a contraction whenever $\sigma=\Re(s)>\delta$.
Notice also that using the variational principle for the topological pressure, it is possible to show that there exist $a_0,b_0>0$
such that for all $\sigma \in \R$,
\begin{equation}
\label{Pressure2}
P(\sigma)\leq a_0-\sigma b_0.
\end{equation}

We continue with a key determinantal identity. We point out that representations of Selberg zeta functions
as Fredholm determinants of transfer operators have a long history going back to Fried \cite{Fried1}, Pollicott 
\cite{Pollicott1} and also Mayer \cite{Mayer1,Mayer2} for the Modular surface. For more recent works involving transfer operators and unitary representations we also mention \cite{Pohl1,Pohl2}.
\begin{propo}
For all $\Re(s)$ large, we have the identity :
\begin{equation}
\label{eq1}
\det(I-\lt_{\rr,s})=L_\Gamma(s,\rr),
\end{equation}
\end{propo}
\noindent{\it Proof}. Remark that the above statement implies analytic continuation to $\C$ of each L-function $L_\Gamma(s,\rr)$, since each $s\mapsto \det(I-\lt_{\rr,s})$ is readily an entire function of $s$. 
For all integer $N\geq 1$, let us compute the trace of $\lt_{\rr,s}^N$. 
Our basic reference for the theory of Fredholm determinants on Hilbert spaces is \cite{Simon}. Let $(e_1,\ldots,e_{d_\rr})$ be an orthonormal
basis of $V_\rr$. For each disc $\D_j$ let $(\varphi_\ell^{j})_{\ell \in \N}$ be a Hilbert basis of the Bergmann space $H^2(\D_j)$, that
is the space of square integrable holomorphic functions on $\D_j$. Then the family defined by
$$\Psi_{j,\ell,k}(z):=\left \{ \varphi_\ell^{j}(z)e_k\ \mathrm{if}\ z\in \D_j \atop 0\  \mathrm{otherwise}, \right.$$
is a Hilbert basis of $H^2_\rr(\Omega)$. Writing
$$\langle \lt_{\rr,s}^N(\Psi_{j,\ell,k}) ,\Psi_{j,\ell,k}\rangle_{H^2_\rr(\Omega)}=
\sum_{\alpha \in \mathscr{W}_N^j} \int_{\D_j} (\gamma_\alpha'(z))^s \varphi^j_\ell(\gamma_\alpha z)\overline{\varphi^j_\ell(z)}
\langle e_k \rr(\gamma_\alpha), e_k \rangle_\rr dm(z),$$
we deduce that
$$\mathrm{Tr}( \lt_{\rr,s}^N)=\sum_{j,\ell,k}\langle \lt_{\rr,s}^N(\Psi_{j,\ell,k}) ,\Psi_{j,\ell,k}\rangle_{H^2_\rr(\Omega)}$$
$$=\sum_j \sum_{\alpha \in \mathscr{W}_N^j \atop \alpha_1=m+j} \chi_\rr(\gamma_\alpha) \int_{\D_j} (\gamma_\alpha'(z))^s B_{\D_j}(\gamma_\alpha z,z)dm(z),$$
where $\chi_\rr$ is the character of $\rr$ and $B_{\D_j}(w,z)$ is the {\it Bergmann reproducing kernel} of $H^2(\D_j)$.
There is an explicit formula for the Bergmann kernel of a disc $\D_j=D(c_j,r_j)$ : 
$$B_{\D_\ell}(w,z)=\frac{r_j^2}{\pi \left [ r_j^2-(w-c_j)(\overline{z}-c_j)\right ]^2}.$$
It is now an exercise involving Stoke's and Cauchy formula (for details we refer to Borthwick \cite{Borthwick}, P. 306) to obtain the Lefschetz identity
$$\int_{\D_j} (\gamma_\alpha'(z))^s
 B_{\D_j}(\gamma_\alpha z,z) dm(z)= \frac{(\gamma_\alpha'(x_\alpha))^s}{1-\gamma_\alpha'(x_\alpha)},$$
 where $x_\alpha$ is the unique fixed point of $\gamma_\alpha:\D_j\rightarrow \D_j$. Moreover, 
 $$\gamma_\alpha'(x_\alpha)=e^{-l( \mathcal{C}_\alpha)},$$
 where $\mathcal{C}_\alpha$ is the closed geodesic represented by the conjugacy class of $\gamma_\alpha \in \Gamma$, and
 $l( \mathcal{C}_\alpha)$ is the length. There is a one-to-one correspondence between prime reduced words (up to circular permutations) in
 $$\bigcup_{N\geq 1} \bigcup_{j=1}^{2m} \{\alpha \in \mathscr{W}_N^j\ \mathrm{such\ that}\ \alpha_1=m+j \},$$
and prime conjugacy classes in $\Gamma$ (see Borthwick \cite{Borthwick}, P. 303), therefore each prime conjugacy class in $\Gamma$ and its iterates appear in the above sum, when $N$ ranges from $1$ to $+\infty$.
 
We have therefore reached formally (absolute convergence is valid for $\Re(s)$ large, see later on)
$$\sum_{N\geq 1} \frac{1}{N}\mathrm{Tr}( \lt_{\rr,s}^N)=\sum_{N\geq 1} \frac{1}{N}\sum_j \sum_{\alpha \in \mathscr{W}_N^j \atop \alpha_1=m+j}
\chi_\rr(\gamma_\alpha) \frac{(\gamma_\alpha'(x_\alpha))^s}{1-\gamma_\alpha'(x_\alpha)}$$
$$=\sum_{\mathcal{C}\in \mathcal{P}} \sum_{k\geq 1} \frac{\chi_\rr(\mathcal{C}^k)}{k}
\frac{e^{-skl(\mathcal{C})}}{1-e^{-kl(\mathcal{C})}}.$$
The prime orbit theorem for convex co-compact groups says that as $T\rightarrow +\infty$, (see for example \cite{Lalley,Naud3}), 
$$\#\{ (k,\mathcal{C})\in \N_0\times \mathcal{P}\ :\ kl(\mathcal{C})\leq T\}=\frac{e^{\delta T}}{\delta T}\left (1+o(1)\right).$$
On the other hand, since $\chi_\rr$ takes obviously finitely many values on $\G$
we get absolute convergence of the above series for $\Re(s)>\delta$. For all $\Re(s)$ large, we get again formally
$$\det(I-\lt_{\rr,s})=\mathrm{exp}\left( \sum_{N\geq 1} \frac{1}{N}\mathrm{Tr}( \lt_{\rr,s}^N)   \right)$$
$$=\mathrm{exp}\left(-\sum_{\mathcal{C},k,n} \frac{\chi_\rr(\mathcal{C}^k)}{k}
e^{-(s+n)kl(\mathcal{C})} \right)=\prod_{\mathcal{C}\in \mathcal{P}}\prod_{n\in \N}
\mathrm{exp}\left(-\sum_{k\geq 1} \frac{\chi_\rr(\mathcal{C}^k)}{k} e^{-(s+n)kl(\mathcal{C})} \right)$$
$$=\prod_{\mathcal{C}\in \mathcal{P}} \prod_{k\in \N} \det \left (Id_{V_\rr}-\rr(\mathcal{C^k})e^{-(s+k)l(\mathcal{C})} \right).$$
This formal manipulations are justified for $\Re(s)>\delta$ by using the spectral radius estimate (\ref{radius1}) and the fact that if $A$
is a trace class operator on a Hilbert space $\mathcal{H}$ with $\Vert A\Vert_{\mathcal{H}}<1$ then we have 
$$\det(I-A)=\mathrm{exp}\left( -\sum_{N\geq 1} \frac{1}{N}\mathrm{Tr}( A^N)\right),$$
(this is a direct consequence of Lidskii's theorem, see \cite{Simon}, chapter 3).
The proof is finished and we have claim $1)$ of Theorem \ref{main1}. $\square$

\bigskip
\noindent
Claim $3)$ follows from the formula (valid for $\Re(s)>\delta$)
$$ \det(I-\lt_{\rr,s})=\mathrm{exp}\left(-\sum_{\mathcal{C},k,n} \frac{\chi_\rr(\mathcal{C}^k)}{k}
e^{-(s+n)kl(\mathcal{C})} \right),$$
and the identity for the character of the regular representation (see \cite{serre}, chapter 2) 
\begin{equation}
\label{Schur1}
\sum_{\rr\  \mathrm{irreducible}} d_\rr \chi_\rr(g)=\vert \G \vert \mathcal{D}_e(g),
\end{equation}
where $\mathcal{D}_e$ is the dirac mass at the neutral element $e$. Indeed, using (\ref{Schur1}), we get
$$\prod_{\rr\  \mathrm{irreducible}} \left (\det(I-\lt_{\rr,s})\right)^{d_\rr}=
\mathrm{exp}\left(-\vert \G \vert \sum_{k,n} \sum_{\mathcal{C} \in \mathcal{P}\atop \mathcal{C}\equiv Id\ \mathrm{mod}\ p}\frac{1}{k}
e^{-(s+n)kl(\mathcal{C})} \right).$$
To see that this is exactly the Euler product defining $Z_{\Gamma(p)}(s)$, observe that since for $p$ large we have (by Gamburd's result \cite{Gamburd1})
$$\Gamma\backslash \Gamma(p)\simeq \G,$$
each conjugacy class in $\Gamma$ of elements belonging to $\Gamma(p)$ {\it splits} into $\vert \G \vert$ conjugacy classes in
$\Gamma(p)$. The details of the group theoretic arguments are in \cite{JN1}, section 2, and it rests on the fact that the only abelian subgroups of $\Gamma$ are elementary subgroups.

\subsection{Singular value estimates}
The proof of claim $2)$ will require more work and will use singular values estimates for vector-valued operators.  We now recall a few facts on singular values of trace class operators.
Our reference for that matter is for example the book \cite{Simon}. If $T:\mathcal{H}\rightarrow \mathcal{H}$ is a compact operator acting on a Hilbert space $\mathcal{H}$, the {\it singular value sequence} is by definition the sequence $\mu_1(T)=\Vert T\Vert\geq \mu_2(T)\geq\ldots \geq \mu_n(T)$ of the eigenvalues of the positive self-adjoint operator $\sqrt{T^*T}$. To estimate singular values in a vector valued setting, we will rely on the following fact.
\begin{lem}
\label{singular}
Assume that $(e_j)_{j\in J}$ is a Hilbert basis of $\mathcal{H}$, indexed by a countable set $J$. Let $T$ be a compact operator on $\mathcal{H}$.
Then for all subset $I\subset J$ with $\# I=n$ we have 
$$\mu_{n+1}(T)\leq \sum_{j\in J\setminus I} \Vert T e_j \Vert_{\mathcal H}.$$
\end{lem}
\noindent {\it Proof}. By the min-max principle for bounded self-adjoint operators, we have
$$\mu_{n+1}(T)=\min_{\mathrm{dim}(F)=n} \max_{w \in F^{\perp},\Vert w \Vert=1} \langle \sqrt{T^*T}w,w\rangle.$$
Set $F=\mathrm{Span}\{ e_j,\ j\in I\}$. Given $w=\sum_{j\not \in I} c_j e_j$ with $\sum_j \vert c_j \vert^2=1$, we obtain via Cauchy-Schwarz inequality
$$\vert \langle \sqrt{T^*T}w,w\rangle\vert \leq \Vert \sqrt{T^*T}(w)\Vert=\Vert T(w)\Vert\leq \sum_{j\not \in I} \Vert T(e_j)\Vert,$$
which concludes the proof. $\square$ 
 
Our aim is now to prove the following bound.
\begin{propo}
\label{eigen1}
Let $(\lambda_k(\lt_{\rr,s}))_{k\geq 1}$ denote the eigenvalue sequence of the compact operators $\lt_{\rr,s}$. There exists $C>0$ and
$0<\eta$ such that for all $s\in \C$ and all representation $\rr$, we have for all $k$,
$$\vert\lambda_k( \lt_{\rr,s}) \vert \leq Cd_\rr e^{C\vert s\vert} e^{-\frac{\eta}{d_\rr}k}.$$
\end{propo}
\noindent Before we prove this bound, let us show quickly how the combination of the above bound with (\ref{radius1}) 
gives the estimate $2)$ of Theorem
\ref{main1}. By definition of Fredholm determinants, we have
$$\log \vert L_\Gamma(s,\rr) \vert\leq \sum_{k=1}^\infty \log(1+\vert\lambda_k( \lt_{\rr,s}) \vert)$$
$$=\sum_{k=1}^N \log(1+\vert\lambda_k( \lt_{\rr,s}) \vert) +\sum_{k=N+1}^\infty \log(1+\vert\lambda_k( \lt_{\rr,s}) \vert),$$
where $N$ will be adjusted later on. The first term is estimated via (\ref{Pressure2}) as
$$ \sum_{k=1}^N \log(1+\vert\lambda_k( \lt_{\rr,s}) \vert)\leq \widetilde{C}(\vert s\vert +1)N,$$
for some large constant $\widetilde{C}>0$. On the other hand we have by the eigenvalue bound from Proposition \ref{eigen1}
$$\sum_{k=N+1}^\infty \log(1+\vert\lambda_k( \lt_{\rr,s}) \vert )\leq \sum_{k=N+1}^\infty \vert\lambda_k( \lt_{\rr,s}) \vert$$
$$\leq Cd_\rr e^{C\vert s\vert} \sum_{k\geq N+1} e^{-\frac{\eta}{d_\rr}k}=
Cd_\rr e^{C\vert s\vert} \frac{e^{-(N+1)\eta/d_\rr}}{1-e^{-\eta/d_\rr}}$$
$$\leq C' \frac{d_\rr^2}{\eta}e^{C\vert s\vert}  e^{-N\frac{\eta}{d_\rr}}.$$
Choosing $N=B[\vert s\vert d_\rr]+B[d\rr \log d_\rr]$ for some large $B>0$ leads to
$$\sum_{k=N+1}^\infty \log(1+\vert\lambda_k( \lt_{\rr,s}) \vert )\leq \widetilde{B}$$
for some constant $\widetilde{B}>0$ uniform in $\vert s\vert$ and $d_\rr$. Therefore we get 
$$ \log \vert L_\Gamma(s,\rr) \vert\leq O\left(d_\rr \log(d_\rr)(\vert s\vert^2+1)\right),$$
which is the bound claimed in statement $2)$.

\noindent {\it Proof of Proposition \ref{eigen1}}. We first recall that if $\D_j=D(c_j,r_j)$, an explicit Hilbert basis of the Bergmann space $H^2(\D_j)$ is given by
the functions ( $\ell=0,\ldots,+\infty$, $j=1,\ldots,2m$)
$$\varphi_\ell^{(j)}(z)=\sqrt{\frac{\ell+1}{\pi}}\frac{1}{r_j} \left (\frac{z-c_j}{r_j} \right)^\ell.$$
By the Schottky property, one can find $\eta_0>0$ such for all $z\in \D_j$, for all $i\neq j$ we have $\gamma_i(z)\in \D_{i+m}$ and
$$\frac{\vert \gamma_i(z)-c_{m+i}\vert}{r_{m+i}}\leq e^{-\eta_0},$$
so that we have uniformly in $i,z$,
\begin{equation}
\label{Ucont2}
 \vert \varphi_\ell^{(i+m)}(\gamma_i z)\vert \leq Ce^{-\eta_1 \ell},
\end{equation}
for some $0<\eta_1<\eta_0$. Going back to the basis $\Psi_{j,\ell,k}(z)$ of $H^2_\rr(\Omega)$, we can write
$$\Vert \lt_{\rr,s}( \Psi_{j,\ell,k})\Vert^2_{H^2_\rr}=\sum_{n=1}^{2m} \sum_{i,i'\neq n}
\int_{\D_n}(\gamma_i(z))^s \overline{(\gamma_{i'}(z))^s} 
\langle\Psi_{j,\ell,k}(\gamma_i z)\rr(\gamma_i),\Psi_{j,\ell,k}(\gamma_{i'} z)\rr(\gamma_{i'})  \rangle_\rr dm(z).$$
Using Schwarz inequality and unitarity of the representation $\rr$ for the inner product $\langle .,.\rangle_\rr$, 
we get by (\ref{Ucont2}) and also (\ref{Bound1}),
$$ \Vert \lt_{\rr,s}( \Psi_{j,\ell,k})\Vert^2_{H^2_\rr}\leq \widetilde{C}e^{\widetilde{C}\vert s \vert} e^{-2\eta_1 \ell},$$
for some large constant $\widetilde{C}>0$. We can now use Lemma \ref{singular} to write
$$\mu_{2md_\rho n+1}(\lt_{\rr,s})\leq \sum_{j=1}^{2m} \sum_{\ell=n}^{+\infty} \sum_{k=1}^{d_\rr} \Vert \lt_{\rr,s}(\Psi_{j,\ell,k}) \Vert_{H^2_\rr} $$
$$\leq C d_\rho e^{\widetilde{C}\vert s\vert} e^{-\eta_1 n},$$
for some $C>0$. Given $N\in \N$, we write $N=2md_\rr k+r$ where $0\leq r<2md_\rr$ and $k=[\frac{N}{2md_\rho}]$. We end up with
$$\mu_{N+1}(\lt_{\rho,s})\leq\mu_{2md_\rr k+1}(\lt_{\rr,s})\leq C'd_\rr e^{\widetilde{C}\vert s\vert} e^{-\eta_2 N/d_\rr},$$
for some $\eta_2>0$. To produce a bound on the eigenvalues, we use then a variant of Weyl inequalities (see \cite{Simon}, Thm 1.14) to get
$$\vert \lambda_N(\lt_{\rr,s})\vert\leq \prod_{k=1}^N \vert \lambda_k(\lt_{\rr,s})\vert \leq \prod_{k=1}^N \mu_k(\lt_{\rr,s}),$$
which yields
$$\vert \lambda_N(\lt_{\rr,s})\vert \leq C_1 d_\rr e^{C_2 \vert s\vert} e^{-\frac{\eta_2}{Nd_\rr}\sum_{k=1}^N k}.$$
Using the well known identity $\sum_{k=1}^N k=\frac{N(N+1)}{2}$ we finally recover
$$\vert \lambda_N(\lt_{\rr,s})\vert \leq C_1 d_\rr e^{C_2 \vert s\vert} e^{-\frac{\eta N}{d_\rr}},$$
for some $\eta>0$ and the proof is done. $\square$

The reader will have noticed that we have used no specific information at all about the representation theory of $\G=SL_2(\Fp)$ and that the above routine works verbatim for abstract groups $\G$ as long as we have a natural group homomorphism $\Gamma \rightarrow \G$. 
\section{Zero-free regions for $L$-functions and explicit formulae}
The goal of this section is to prove the following result which will allows us to convert zero-free regions into upper bounds on sums
over closed geodesics.

\begin{propo}
\label{ExplicitF}
Fix $\overline{\alpha}>0$, $0\leq \sigma<\delta$ and $\varepsilon>0$. Then there exists a $C_0^\infty$ test function $\varphi_0$, with
$\varphi_0\geq 0$, $\mathrm{Supp}(\varphi_0)=[-1,+1]$ and such that that for $\rr$ non trivial, if $L_\Gamma(s,\rr)$ has no zeros in the rectangle
$$\{ \sigma \leq \Re(s) \leq 1\ \mathrm{and}\ \vert \Im(s)\vert \leq (\log T)^{1+\overline{\alpha}}\},$$
for some $T$ large enough, then we have 
$$\sum_{\mathcal{C},k} \chi_\rr(\mathcal{C}^k)
\frac{l(\mathcal{C})}{1-e^{kl(\mathcal{C})}}\varphi_0\left( \frac{kl(\mathcal{C})}{T}\right)=
O\left(d_\rr \log(d_\rr)e^{(\sigma+\varepsilon)T}\right),$$
where the implied constant is uniform in $T, d_\rr$.
\end{propo}
The proof will occupy the full section and will be broken into several elementary steps.

\bigskip
\subsection{Preliminary Lemmas}
We start this section by the following fact from harmonic analysis.
\begin{lem}
\label{Fourier1}
For all $\alpha>0$, there exists $C_1,C_2>0$ and a positive test function $\varphi_0 \in C_0^\infty(\R)$ with 
$\mathrm{Supp}(\varphi)=[-1,+1]$ such that for all $\vert \xi\vert \geq 2$, we have
$$\vert \widehat{\varphi_0}(\xi)\vert \leq C_1 e^{\vert \Im(\xi)\vert} \exp\left (-C_2\frac{\vert \Re(\xi)\vert}{(\log\vert \Re(\xi)\vert)^{1+\alpha}}   \right),$$
where $\widehat{\varphi_0}(\xi)$ is the Fourier transform, defined as usual by
$$ \widehat{\varphi_0}(\xi)=\int_{-\infty}^{+\infty} \varphi_0(x)e^{-ix\xi}dx. $$
\end{lem}
\noindent {\it Proof.} It is known from the Beurling-Malliavin multiplier Theorem, or the Denjoy-Carleman Theorem, that for
compactly supported test functions $\psi$, one cannot beat the rate ($\xi \in \R$, large)
$$ \vert \widehat{\psi}(\xi)\vert =O\left (\exp\left(-C\frac{\vert \xi\vert}{\log\vert \xi\vert}\right)  \right),$$
because this rate of Fourier decay implies quasi-analyticity (hence no compactly supported test functions). We refer the reader to
\cite{Katz}, chapter 5 for more details. The above statement is definitely a {\it folklore} result. However since we need a precise control for complex valued
$\xi$ and couldn't find the exact reference for it, we provide an outline of the proof which follows closely the construction that one can find in \cite{Katz}, chapter 5, Lemma 2.7.

Let $(\mu_j)_{j\geq 1}$ be a sequence of positive numbers such that $\sum_{j=1}^\infty \mu_j=1$. For all $k\in \Z$, set
$$\varphi_N(k)=\prod_{j=1}^N \frac{\sin(\mu_j k)}{\mu_j k},\ \  \varphi(k)=\prod_{j=1}^\infty \frac{\sin(\mu_j k)}{\mu_j k}.$$
Consider the Fourier series given by
$$f(x):=\sum_{k \in \Z} \varphi(k) e^{ikx},\ \ f_N(x):=\sum_{k \in \Z} \varphi_N(k) e^{ikx},$$
then one can observe that by rapid decay of $\varphi(k)$, $f(x)$ defines a $C^\infty$ function on $[-2\pi,2\pi]$.
On the other hand, one can check that $f_N(x)$ converges uniformly to $f$ as $N$ goes to $\infty$ and that
$$f_N(x)=(g_1 \star g_2\star \ldots \star g_N)(x),$$
where $\star$ is the convolution product and each $g_j$ is given by
$$g_j(x):=\left \{ \frac{2\pi}{\mu_j}\ \mathrm{if}\ \vert x\vert \leq \mu_j \atop 0\ \mathrm{elsewhere}. \right.$$
From this observation one deduces that $f$ is positive and supported on $[-1,+1]$ since we assume
$$\sum_{j=1}^\infty \mu_j=1.$$
We now extend $f$ outside $[-1,+1]$ by zero and write by integration by parts and Schwarz inequality,
$$\vert \widehat{ f}(\xi)\vert \leq  \frac{e^{\vert \Im(\xi)\vert}}{\vert \Re(\xi)\vert^N}\Vert f^{(N)}\Vert_{L^2(-1,+1)}.$$
By Plancherel formula, we get
$$\Vert f^{(N)}\Vert_{L^2(-1,+1)}^2=\sum_{k \in \Z} k^{2N}(\varphi(k))^2\leq C \prod_{j=1}^{N+1}\mu_j^{-2},$$
where $C>0$ is some universal constant. Fixing $\epsilon>0$, we now choose
$$\mu_j=\frac{\widetilde{C}}{j(\log(1+j))^{1+\epsilon}},$$
where $\widetilde{C}$ is adjusted so that $\sum_{j=1}^\infty \mu_j=1$, and we get
$$\vert \widehat{ f}(\xi)\vert \leq  \frac{e^{\vert \Im(\xi)\vert}}{\vert \Re(\xi)\vert^N} (C_1)^N N! e^{N(1+\epsilon)\log\log(N)}.$$
Using Stirling's formula and choosing $N$ of size
$$N=\left [ \frac{\vert \Re(\xi)\vert}{(\log(\vert\Re(\xi)\vert)^{1+2\epsilon}} \right ] $$
yields (after some calculations) to
$$\vert \widehat{ f}(\xi)\vert \leq  O\left (e^{\vert \Im(\xi)\vert} e^{-C_2 \frac{\vert \Re(\xi)\vert}{(\log(\vert \Re(\xi)\vert)^{1+2\epsilon}}}  
\right),$$
and the proof is finished. $\square$

One can obviously push the above construction further below the threshold $\frac{\vert \xi \vert}{\log\vert \xi\vert}$
by obtaining decay rates of the type
$$\exp \left (-\frac{\vert \xi\vert}{\log\vert \xi \vert \log(\log\vert \xi \vert) \ldots (\log_{(n)}\vert \xi \vert)^{1+\alpha}}\right),$$
where $\log_{(n)}(x)=\log\log\ldots\log(x)$, iterated $n$ times. However this would only yield a very mild improvement to the main statement,
so we will content ourselves with the above lemma.

We continue with another result which will allow us to estimate the size of the log-derivative of $L_\Gamma(s,\rr)$ in a narrow rectangular
zero-free region. More precisely, we have the following:
\begin{propo}
\label{Derivative1}
Fix $\sigma<\delta$. For all $\epsilon>0$, there exist $C(\epsilon), R(\epsilon) >0$ such that for all $R\geq R(\epsilon)$,
if $L_\Gamma(s,\rr)$ ($\rr$ is non trivial) has no zeros in the rectangle
$$\{ \sigma \leq \Re(s) \leq 1\ \mathrm{and}\ \vert \Im(s)\vert \leq R\},$$
then we have for all $s$ in the smaller rectangle
$$\{ \sigma+\epsilon \leq \Re(s) \leq 1\ \mathrm{and}\ \vert \Im(s)\vert \leq C(\epsilon)R\},$$
$$\left \vert \frac{L'_\Gamma(s,\rr)}{L_\Gamma(s,\rr)}   \right \vert \leq B(\epsilon) d_\rr \log(d_\rr) R^6.$$
\end{propo}
\noindent {\it Proof}. We will use Caratheodory's Lemma and take advantage of the {\it a priori} bound from Theorem \ref{main1}. More precisely,
our goal is to rely on this estimate (see Titchmarsh \cite{Tits}, 5.51).
\begin{lem}
\label{Bigtits}
 Assume that $f$ is a holomorphic function on a neighborhood of the closed disc $\overline {D}(0,r)$, then
 for all $r'<r$, we have
 $$ \max_{\vert z\vert\leq r'} \vert f'(z) \vert \leq \frac{8r}{(r-r')^2}\left ( 
 \max_{\vert z  \vert \leq r} \vert \Re(f(z))\vert+\vert f(0)\vert \right).$$
 \end{lem}
 First we recall that for all $\Re(s)>\delta$, $L_\Gamma(s,\rr)$ does not vanish and has a representation as
 $$L_\Gamma(s,\rr)=\mathrm{exp}\left(-\sum_{\mathcal{C},k} \frac{\chi_\rr(\mathcal{C}^k)}{k}
\frac{e^{-skl(\mathcal{C})}}{1-e^{kl(\mathcal{C})}} \right),$$
so that we get for all $\Re(s)\geq A>\delta$,
\begin{equation}
\label{Center1}
\left \vert \log\vert L_\Gamma(s,\rr)\vert \right \vert \leq C_A d_\rr,\ \ \left \vert \frac{L'_\Gamma(s,\rr)}{L_\Gamma(s,\rr)} \right \vert \leq  C'_A d_\rr
\end{equation}
where $C_A, C'_A>0$ are uniform constants on all half-planes 
$$\{\Re(s)\geq A>\delta \}.$$
We have simply used the prime orbit theorem and the trivial bound on characters of unitary representations:
 $$\vert\chi_\rr(g)\vert \leq d_\rr,\ \mathrm{for\ all}\  g\in \G.$$
 Let us now assume that $L_\Gamma(s,\rr)$ does not vanish on the rectangle
 $$\{ \sigma \leq \Re(s) \leq 1\ \mathrm{and}\ \vert \Im(s)\vert \leq R\}.$$ 
 Consider the disc $D(M,r)$ centered at $M$ and with radius $r$ where $M(\sigma,R)$ and $r(\sigma,R)$ are given by
 $$M(\sigma,R)=\frac{R^2}{2(1-\sigma)}+\frac{\sigma+1}{2};\ \ r(\sigma,R)=M(\sigma,R)-\sigma,$$
 see the figure below.
 
 \begin{center}
\includegraphics[scale=0.5]{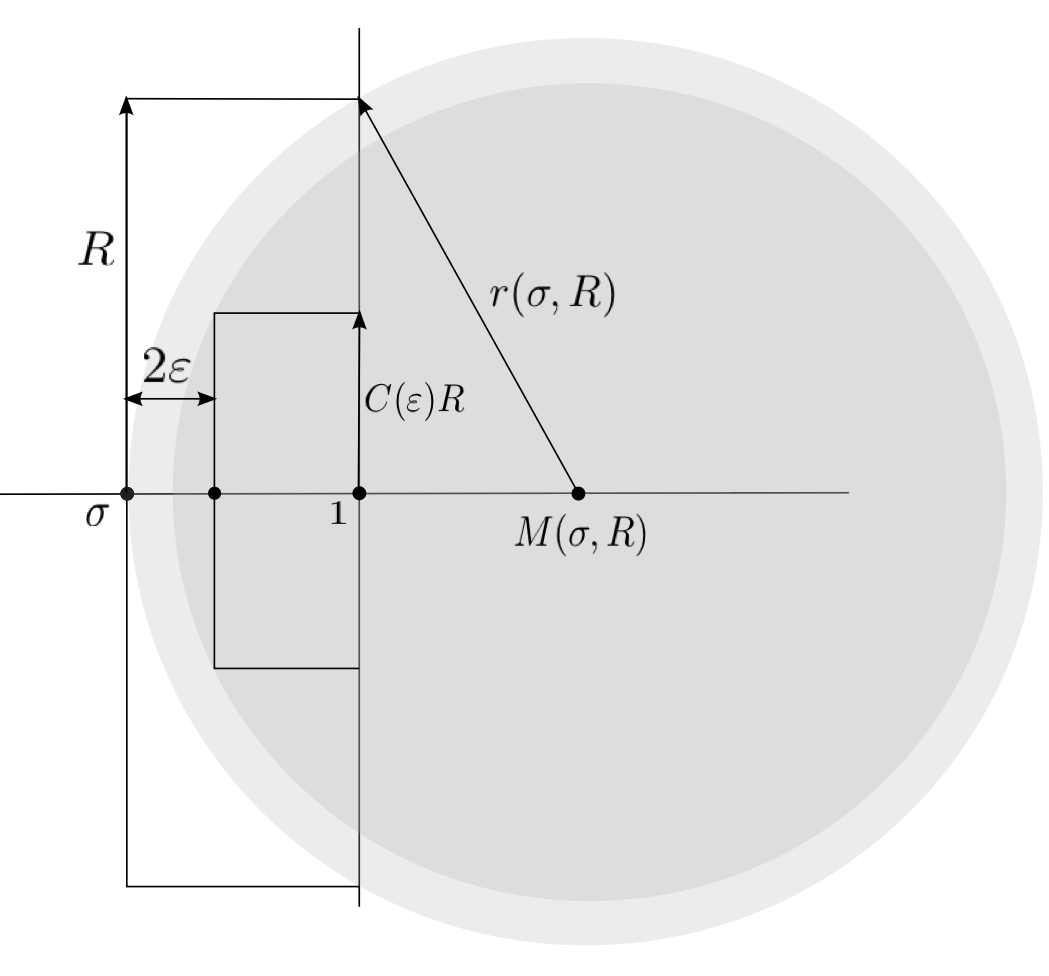}
\end{center}
Since by assumption $s\mapsto L_\Gamma(s,\rr)$ does not vanish on the closed disc $\overline{D}(M,r)$, we can choose a determination
of the complex logarithm of $L_\Gamma(s,\rr)$ on this disc to which we can apply Lemma \ref{Bigtits} on the smaller disc 
$D(M,r-\varepsilon)$, which yields (using the a priori bound from Theorem \ref{main1} and estimate (\ref{Center1}))
$$\left \vert \frac{L'_\Gamma(s,\rr)}{L_\Gamma(s,\rr)}   \right \vert
 \leq C\frac{r}{\varepsilon}\left (d_\rr \log(d_\rr) r^2 +A_1d_\rr \right )$$
 $$=O\left ( R^6 d_\rr \log(d_\rr) \right),$$
 where the implied constant is uniform with respect to $R$ and $d_\rr$.
 Looking at the picture, the smaller disc $D(M,r-\varepsilon)$ contains a rectangle
 $$\{ \sigma+2\varepsilon \leq \Re(s) \leq 1\ \mathrm{and}\ \vert \Im(s)\vert \leq L(\varepsilon)\},$$
 where $L(\varepsilon)$ satisfies the identity (Pythagoras Theorem!)
 $$L^2(\varepsilon)=\varepsilon(2M-2\sigma-3\varepsilon),$$
which shows that 
$$ L(\epsilon)\geq C(\varepsilon)R,$$
with $C(\varepsilon)>0$, as long as $R\geq R_0(\epsilon)$, for some $R_0>0$. The proof is done. $\square$

\subsection{Proof of the Proposition \ref{ExplicitF}}
We are now ready to prove the main result of this section, by combining the above facts with a standard contour deformation argument. We fix a small $\varepsilon>0$ and $0<\alpha<\overline{\alpha}$. We use Lemma
\ref{Fourier1} to pick a test function $\varphi_0$ with Fourier decay as described, with same exponent $\alpha$. We set for all $T>0$, and $s\in \C$,
$$\psi_T(s)=\int_{-\infty}^{+\infty} e^{sx}\varphi_0\left(\frac{x}{T}\right)dx$$
$$=T\widehat{\varphi_0}(isT).$$
By the estimate from Lemma \ref{Fourier1}, 
we have 
\begin{equation}
\label{testdecay}
\vert \psi_T(s)\vert\leq C_1 Te^{T\vert \Re(s)\vert}\exp\left(-C_2 \frac{\vert \Im(s)\vert T}{(\log(T\vert \Im(s)\vert)^{1+\alpha}}\right).
\end{equation} 
We fix now $A>\delta$ and consider the contour integral
$$I(\rr,T)=\frac{1}{2i\pi} \int_{A-i\infty}^{A+i\infty} \frac{L'_\Gamma(s,\rr)}{L_\Gamma(s,\rr)}\psi_T(s)ds.$$
Convergence is guaranteed by estimate (\ref{Center1}) and rapid decay of $\vert \psi_T(s)\vert$ on vertical lines.
Because we choose $A>\delta$, we have absolute convergence of the series
$$\frac{L'_\Gamma(s,\rr)}{L_\Gamma(s,\rr)}=\sum_{\mathcal{C},k} \chi_\rr(\mathcal{C}^k)
\frac{l(\mathcal{C})e^{-skl(\mathcal{C})}}{1-e^{kl(\mathcal{C})}}$$
on the vertical line $\{\Re(s)=A\}$, 
and we can use Fubini to write
$$I(\rr,T)=\sum_{\mathcal{C},k} \chi_\rr(\mathcal{C}^k)
\frac{l(\mathcal{C})}{1-e^{kl(\mathcal{C})}}e^{-Akl(\mathcal{C})}\frac{1}{2\pi}\int_{-\infty}^{+\infty} e^{-itkl(\mathcal{C})}
\widehat{\psi_T}(iA-t)dt,$$
and Fourier inversion formula gives
$$I(\rr,T)=\sum_{\mathcal{C},k} \chi_\rr(\mathcal{C}^k)
\frac{l(\mathcal{C})}{1-e^{kl(\mathcal{C})}}\varphi_0\left( \frac{kl(\mathcal{C})}{T}\right).$$
Assuming that $L_\Gamma(s,\rr)$ has no zeros in 
$$\{ \sigma \leq \Re(s) \leq 1\ \mathrm{and}\ \vert \Im(s)\vert \leq R\},$$
where $R$ will be adjusted later on, our aim is to use Proposition \ref{Derivative1} to deform the contour integral $I(\rr,T)$ as depicted in the figure
below. 
 \begin{center}
\includegraphics[scale=0.6]{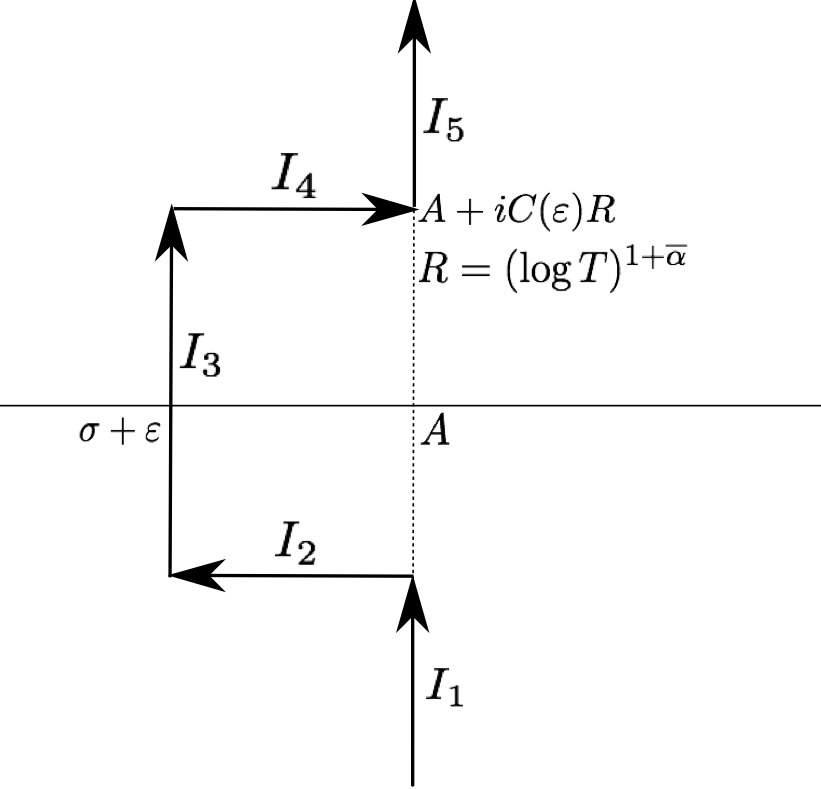}
\end{center} 
 Writing $I(\rr,T)=\sum_{j=1}^5 I_j$ (see the above figure), we need to estimate carefully each contribution. In the course of the proof, we will use the following basic fact.
 \begin{lem}
 \label{expo1}
 Let $\phi:[M_0,+\infty)\rightarrow \R^+$ be a $C^2$ map with $\phi'(x)>0$ on $[M_0,+\infty)$ and satisfying
 $$(*)\ \ \ \ \sup_{x\geq M_0} \left \vert \frac{\phi''(x)}{(\phi'(x))^3}\right \vert \leq C, $$
 then we have for all $M\geq M_0$,
 $$\int_M^{+\infty} e^{-\phi(t)}dt\leq \frac{e^{-\phi(M)}}{\phi'(M)}+Ce^{-\phi(M)}.$$
 \end{lem}
 {\noindent \it Proof.} First observe that condition $(*)$ implies that 
 $$x\mapsto \frac{1}{(\phi'(x))^2}$$
 has a uniformly bounded derivative, which is enough to guarantee that 
 $$\lim_{x\rightarrow +\infty} \frac{e^{-\phi(x)}}{\phi'(x)}=0. $$
 In particular $\lim_{x\to +\infty} \phi(x)=+\infty$ and for all $M\geq M_0$, $\phi:[M,+\infty)\rightarrow [\phi(M),+\infty)$ is a $C^2$-diffeomorphism.
 A change of variable gives
 $$\int_M^{+\infty} e^{-\phi(t)}dt=\int_{\phi(M)}^{+\infty} e^{-u}\frac{du}{\phi'(\phi^{-1}(u))},$$
 and integrating by parts yields the result. $\square$
 
 \bigskip
\begin{itemize}
\item First we start with $I_1$ and $I_5$. Using estimate (\ref{Center1}) combined with (\ref{testdecay}), we have
$$\vert I_5 \vert \leq C d_\rr T e^{TA} \int_{C(\varepsilon)R} ^{+\infty} e^{-C_2 \frac{tT}{(\log(tT))^{1+\alpha}}}dt,$$
which by a change of variable leaves us with
$$\vert I_5 \vert \leq C d_\rr e^{TA} \int_{C(\varepsilon)RT}^{+\infty} e^{-C_2 \frac{u}{(\log(u))^{1+\alpha}}}du.$$
This where we use Lemma \ref{expo1} with 
$$\phi(x)=C_2 \frac{x}{(\log(x))^{1+\alpha}}.$$
Computing the first two derivatives, we can check that condition $(*)$ is fulfilled and therefore
$$\int_{M}^{+\infty} e^{-C_2 \frac{u}{(\log(u))^{1+\alpha}}}\leq C (\log(M))^{1+\alpha} e^{-C_2 \frac{M}{(\log(M))^{1+\alpha}}},$$
for some universal constant $C>0$. We have finally obtained
$$ \vert I_5 \vert \leq C d_\rr e^{TA}(\log(RT))^{1+\alpha} e^{-C_2 \frac{RT}{(\log(RT))^{1+\alpha}}}.$$
Choosing $R=(\log(T))^{1+\overline{\alpha}}$, with $\overline{\alpha}>\alpha$ gives 

$$\vert I_5 \vert =O\left(d_\rr e^{TA}(\log(T))^{1+\alpha} e^{-C_2 T(\log(T))^{\overline{\alpha}-\alpha}}\right)$$
$$=O\left ( d_\rho e^{-BT}\right),$$
where $B>0$ can be taken as large as we want. The exact same estimate is valid for $I_1$.
\item The case of $I_4$ and $I_2$. Here we use the bound from Proposition \ref{Derivative1} and again (\ref{testdecay}) to get
$$ \vert I_4\vert+\vert I_2\vert=O\left ( d_\rr \log(d_\rr) e^{-BT} \right),$$
where $B$ can be taken again as large as we want.
\item We are left with $I_3$ where 
$$I_3=\frac{1}{2\pi}\int_{-C(\epsilon)R}^{+C(\epsilon)R} \frac{L'_\Gamma(\sigma+\varepsilon+it,\rr)}{L_\Gamma(\sigma+\varepsilon+it,\rr)}
\psi_T(\sigma+\varepsilon+it)dt.$$
Using Proposition \ref{Derivative1} and (\ref{testdecay}) we get
$$\vert I_3 \vert =O\left( d_\rr \log(d_\rr) (\log(T))^{7(1+\overline{\alpha})} e^{(\sigma+\epsilon)T}\right).$$
\end{itemize}
Clearly the leading term in the contour integral is provided by $I_3$, and the proof of Proposition \ref{ExplicitF} is now complete. 

\bigskip
We conclude this section by a final observation. If $\rr=id$ is the {\it trivial representation}, then $L_\Gamma(s,id)=Z_\Gamma(s)$ has a zero
at $s=\delta$, thus the best estimate for the contour integral $I(id,T)$ is given by (\ref{Center1}) and (\ref{testdecay}) which yields (by a change of variable)
$$\vert I(id,T) \vert \leq C_Ad_\rr \int_{-\infty}^{+\infty} \vert \psi_T(A+it)\vert dt \leq \widetilde{C}_A d_\rr Te^{TA}
\int_{-\infty}^{+\infty}\exp\left(-C_2 \frac{\vert  t\vert T}{\left(\log(T\vert t\vert)\right)^{1+\alpha}}\right)dt$$
$$=O\left(d_\rr e^{TA}\right).$$
Since $d_\rr=1$ and $A$ can be taken as close to $\delta$ as we want, the contribution from the trivial representation is of size
\begin{equation}
\label{trivial}
I(id,T)=O\left(e^{(\delta+\epsilon)T} \right).
\end{equation}

\section{Existence of "low lying" zeros for $L_\Gamma(s,\rr)$}
\subsection{Conjugacy classes in $\G$.}
In this section, we will use more precise knowledge on the group structure of 
$$\G=SL_2(\Fp).$$
Our basic reference is the book \cite{Suzuki}, see chapter 3, $\S 6$ for much more general statements over finite fields.
We start by describing the conjugacy classes in $\G$. Since we are only interested in the large $p$ behaviour, we will assume that
$p$ is an odd prime strictly bigger than $3$. Conjugacy classes of elements $g\in \G$ are essentially determined by the roots of the characteristic polynomial
$$\det(xI_2-g)=x^2-\mathrm{tr}(g)x+1,$$
which are denoted by $\lambda, \lambda^{-1}$, where $\lambda \in \Fp^{\times}$. There are three different possibilities.
\begin{itemize}
\item $\lambda\neq \lambda^{-1} \in  \Fp^{\times}$. In that case $g$ is diagonalizable over $\Fp$ and $g$ is conjugate to the matrix
$$D(\lambda)=\left (\begin{array}{cc}
\lambda &0\\
0&\lambda^{-1}
\end{array} \right ).$$
The centralizer $\mathcal{Z}(D(\lambda))=\{ h\in \G\ :\ hD(\lambda)h^{-1}=D(\lambda)\}$ is then equal to the "maximal torus"
$$A=\left\{  \left (\begin{array}{cc}
a &0\\
0&a^{-1}
\end{array} \right )\ :\ a \in   \Fp^{\times}\right\},$$
and we have $\vert A\vert=p-1$, the conjugacy class of $g$ has $p(p+1)$ elements.
\item $\lambda\neq \lambda^{-1} \not \in  \Fp^{\times}$. In that case $\lambda$ belongs to $\mathcal{F}\simeq \mathbb{F}_{p^2}$ the unique quadratic extension of $\Fp$. The root $\lambda$ can be written as 
$$ \lambda=a+b\sqrt{\epsilon},\  \lambda^{-1}=a-b\sqrt{\epsilon},$$
where $\{1,\sqrt{\epsilon}\}$ is a fixed $\Fp$-basis of $\mathcal{F}$.
Therefore $g$ is conjugate to 
$$\left (\begin{array}{cc}
a &\epsilon b\\
b& a
\end{array} \right),$$
and $\vert \mathcal{Z}(g)\vert=p+1$, its conjugacy class has $p(p-1)$ elements.
\item $\lambda=\lambda^{-1}\in \{\pm 1\}$. In that case $g$ is non-diagonalizable unless $g\in \mathcal{Z}(\G)=\{\pm I_2\}$, and is conjugate to
$\pm u$ or $\pm u'$ where 
$$u= \left (\begin{array}{cc}
1 & 1\\
0& 1
\end{array} \right),\ u'=\left (\begin{array}{cc}
1 & \epsilon\\
0& 1
\end{array} \right).$$
The centralizer $\mathcal{Z}(g)$ has cardinality $2p$ and the four conjugacy classes have $p(p+1)$ elements.  
\end{itemize}
Using this knowledge on conjugacy classes, one can construct all irreducible representations and write a character table for $\G$,
but we won't need it. There are two facts that we highlight and will use in the sequel:

\begin{enumerate}
 \item For all $g \in \G$, $\vert \mathcal{Z}(g)\vert \geq p-1$.
 \item For all $\rr$ non-trivial we have $d_\rho\geq \frac{p-1}{2}$.
\end{enumerate}
We will also rely on the very important observation below.
\begin{propo}
\label{Conj1}
Let $\Gamma$ be a convex co-compact subgroup of $SL_2(\Z)$ as above.
Fix $0<\beta<2$, and consider the set $\mathcal{E}_T$ of conjugacy classes $\overline{\gamma}\subset \Gamma\setminus \{Id\}$ such that for all 
$\overline{\gamma}\in \mathcal{E}_T$, we have $\l(\gamma)\leq T:=\beta \log(p)$. Then for all $p$ large and all
$\overline{\gamma_1},\overline{\gamma_2}\in \mathcal{E}_T$, the following are equivalent:
\begin{enumerate}
\item $\mathrm{tr}(\gamma_1)=\mathrm{tr}(\gamma_2)$.
\item $\gamma_1$ and $\gamma_2$ are conjugate in $\G$.
\end{enumerate}
\end{propo}
\noindent{\it Proof}. Clearly $(1)$ implies that $\gamma_1$ and $\gamma_2$ have the same trace modulo $p$. Unless
we are in the cases $\mathrm{tr}(\gamma_1)=\mathrm{tr}(\gamma_2)=\pm 2\ \mathrm{mod}\ p$, we know from the above description of conjugacy classes that
they are determined by the knowledge of the trace.  To eliminate these "parabolic mod $p$" cases, we observe that if $\overline{\gamma}\in \mathcal{E}_T$ satisfies 
$\mathrm{tr}(\gamma)=\pm 2 +kp$ with $k\neq 0$, then 
$$2\cosh(l(\gamma)/2)=\vert \mathrm{tr}(\gamma)\vert\geq p-2,$$
and we get
$$p-2\leq 1+p^{\frac{\beta}{2}},$$
which leads to an obvious contradiction if $p$ is large, therefore $k=0$.  Then it means that 
$\vert\mathrm{tr}(\gamma)\vert=2$ which is impossible since $\Gamma$ has no non trivial parabolic element (convex co-compact hypothesis). Conversely,
if $\gamma_1$ and $\gamma_2$ are conjugate in $\G$, then we have 
$$\mathrm{tr}(\gamma_1)=\mathrm{tr}(\gamma_2)\ \mathrm{mod}\ p.$$
If $\mathrm{tr}(\gamma_1)\neq \mathrm{tr}(\gamma_2)$ then this gives
$$ p\leq \vert \mathrm{tr}(\gamma_1)-\mathrm{tr}(\gamma_2)\vert \leq 4\cosh(T/2)\leq 2(p^{\frac{\beta}{2}}+1),$$
again a contradiction for $p$ large. $\square$
\subsection{Proof of the main result} 
Before we can rigourously prove Theorem \ref{main2}, we need one last fact from representation theory which is a handy folklore formula.
\begin{lem}
\label{Product1}
Let $\G$ be a finite group and let $\rr:\G\rightarrow \mathrm{End}(V_\rr)$ be an irreducible representation.
Then for all $x,y\in \G$, we have
$$\chi_\rr(x)\overline{\chi_\rr(y)}=\frac{d_\rr}{\vert G\vert} \sum_{g\in \G} \chi_\rr(xgy^{-1}g^{-1}).$$
\end{lem}
\noindent {\it Proof.} Writing
$$ \sum_{g\in \G} \chi_\rr(xgy^{-1}g^{-1})=\mathrm{Tr}\left(\rr(x) \sum_{g} \rr(gy^{-1}g^{-1}) \right),$$
we observe that 
$$U_y:=\sum_{g} \rr(gy^{-1}g^{-1})$$
commutes with the irreducible representation $\rr$, therefore by Schur's Lemma \cite{serre} (chapter 2), it has to be of the form
$$U_y=\lambda(y)I_{V_\rr},$$
with $\lambda(y)\in \C$, which shows that
$$\sum_{g\in \G} \chi_\rr(xgy^{-1}g^{-1})=\chi_\rr(x)\lambda(y).$$
Similarly we obtain
$$ \sum_{g\in \G} \chi_\rr(xgy^{-1}g^{-1})=\overline{\chi_\rr(y)}\lambda(x),$$
and evaluating at the neutral element $x=e_\G$ ends the proof since we have
$$U_{e_\G}=\vert \G \vert I_{V_\rr}.\  \square$$

\bigskip We fix some $0\leq \sigma<\delta$. We take $\varepsilon>0$ and $\alpha>0$. We assume that for all non-trivial representation $\rr$, the corresponding $L$-function $L_\Gamma(s,\rr)$ does not vanish on the rectangle 
$$\{ \sigma \leq \Re(s) \leq 1\ \mathrm{and}\ \vert \Im(s)\vert \leq (\log T)^{1+\alpha}\},$$ 
where $T=\beta \log(p)$ with $0<\beta<2$. The idea is to look at the average
$$S(p):=\sum_{\rr\  \mathrm{irreducible}} \vert I(\rr,T) \vert^2,$$
where $I(\rr,T)$ is the sum given by
$$I(\rr,T)=\sum_{\mathcal{C},k} \chi_\rr(\mathcal{C}^k)
\frac{l(\mathcal{C})}{1-e^{kl(\mathcal{C})}}\varphi_0\left( \frac{kl(\mathcal{C})}{T}\right).$$
While each term $I(\rr,T)$ is hard to estimate from below because of the oscillating behaviour of characters, the mean square is tractable thanks to
Lemma \ref{Product1}. Let us compute $S(p)$.
$$S(p)=\sum_{\rr\  \mathrm{irreducible}} \sum_{\mathcal{C},k} \sum_{\mathcal{C'},k'} 
\frac{l(\mathcal{C})l(\mathcal{C'})}{(1-e^{kl(\mathcal{C})})(1-e^{k'l(\mathcal{C'})})}
\varphi_0\left( \frac{kl(\mathcal{C})}{T}\right)\varphi_0\left( \frac{k'l(\mathcal{C'})}{T}\right)
\chi_\rr(\mathcal{C}^k)\overline{\chi_\rr(\mathcal{C'}^{k'})}.$$
Using Lemma \ref{Product1}, we have
$$\chi_\rr(\mathcal{C}^k)\overline{\chi_\rr(\mathcal{C'}^{k'})}=
\frac{d_\rr}{\vert G\vert} \sum_{g\in \G} \chi_\rr(\mathcal{C}^kg(\mathcal{C'})^{-k'}g^{-1}),$$
and Fubini plus the identity
$$\sum_{\rr\  \mathrm{irreducible}} d_\rr \chi_\rr(g)=\vert \G \vert \mathcal{D}_e(g)$$
allow us to obtain
$$S(p)= \sum_{\mathcal{C},k} \sum_{\mathcal{C'},k'} 
\frac{l(\mathcal{C})l(\mathcal{C'})}{(1-e^{kl(\mathcal{C})})(1-e^{k'l(\mathcal{C'})})}
\varphi_0\left( \frac{kl(\mathcal{C})}{T}\right)\varphi_0\left( \frac{k'l(\mathcal{C'})}{T}\right)
\Phi_\G(\mathcal{C}^k,\mathcal{C'}^{k'}),$$
where
$$\Phi_\G(\mathcal{C}^k,\mathcal{C'}^{k'}):=\sum_{g\in \G} \mathcal{D}_e(\mathcal{C}^kg(\mathcal{C'})^{-k'}g^{-1}).$$
Since all terms in this sum are now positive and $\mathrm{Supp}(\varphi_0)=[-1,+1]$, we can fix a small $\varepsilon>0$ and find a constant $C_\varepsilon>0$ such that
$$S(p)\geq C_\varepsilon  \sum_{kl(\mathcal{C})\leq {\tiny T(1-\varepsilon)} \atop k'l(\mathcal{C'})\leq T(1-\varepsilon)} \Phi_\G(\mathcal{C}^k,\mathcal{C'}^{k'}).$$
Observe now that
$$\Phi_\G(\mathcal{C}^k,\mathcal{C'}^{k'})=\sum_{g\in \G} \mathcal{D}_e(\mathcal{C}^k g(\mathcal{C'})^{-k'}g^{-1})\neq 0$$
if and only if $\mathcal{C}^k$ and $\mathcal{C'}^{-k'}$ {\it are in the same conjugacy class mod $p$}, and in that case,
$$\Phi_\G(\mathcal{C}^k,\mathcal{C'}^{k'})=\vert \mathcal{Z}(\mathcal{C}^k)\vert=\vert \mathcal{Z}(\mathcal{C'}^{k'})\vert.$$
Using the lower bound for the cardinality of centralizers, we end up with
$$S(p)\geq C_\varepsilon (p-1) \sum_{\overline{\mathcal{C}^k}=\overline{\mathcal{C'}^{k'}}\ \mathrm{mod}\ p
\atop kl(\mathcal{C}), k'l(\mathcal{C'})\leq {\tiny T(1-\varepsilon)}} 1.$$
Notice that since we have taken $T=\beta\log(p)$ with $\beta<2$, we can use Proposition\ref{Conj1} which says that 
$\mathcal{C}^k$ and $\mathcal{C'}^{-k'}$ are in the same conjugacy class mod $p$ iff they have the same traces (in $SL_2(\Z)$).
It is therefore natural to rewrite the lower bound for $S(p)$ in terms of traces. We need to introduce a bit more notations. Let 
$\mathscr{L}_\Gamma$ be set of traces i.e.
$$\mathscr{L}_\Gamma=\{ \mathrm{tr}(\gamma)\ :\ \gamma \in \Gamma \}\subset \Z.$$
Given $t\in \mathscr{L}_\Gamma$, we denote by $m(t)$ the multiplicity of $t$ in the trace set by
$$m(t)=\#\{ \mathrm{conj\ class\ }\overline{\gamma}\subset\Gamma\ :\ \mathrm{tr}(\gamma)=t \}.$$
We have therefore (notice that multiplicities are squared in the double sum)
$$S(p)\geq C_\varepsilon (p-1) \sum_{t\in \mathscr{L}_\Gamma \atop \vert t \vert\leq {2\cosh(T(1-\varepsilon)/2)} } m^2(t). $$
To estimate from below this sum, we use a trick that goes back to Selberg. By the prime orbit theorem \cite{Naud3, Lalley,Roblin} applied to
the surface $\Gamma \backslash \hh$, we know that for all $T$ large, we have
$$Ce^{(\delta-2\varepsilon)T}\leq \sum_{t\in \mathscr{L}_\Gamma \atop \vert t \vert\leq {2\cosh(T(1-\varepsilon)/2)} } m(t),$$
and by Schwarz inequality we get for $T$ large
$$Ce^{(\delta-2\varepsilon)T}\leq C_0
\left ( \sum_{t\in \mathscr{L}_\Gamma \atop \vert t \vert\leq {2\cosh(T(1-\varepsilon)/2)} } m^2(t)\right)^{1/2} e^{T/4},$$
where we have used the obvious bound
$$\#\{ n\in \Z\ :\ \vert n\vert \leq 2\cosh(T(1-\varepsilon)/2)\}=O(e^{T/2}).$$
This yields the lower bound
$$\sum_{t\in \mathscr{L}_\Gamma \atop \vert t \vert\leq {2\cosh(T(1-\varepsilon)/2)} } m^2(t)\geq C_\varepsilon'  e^{(2\delta-1/2-\varepsilon)T},$$
which shows that one can take advantage of exponential multiplicities in the length spectrum when $\delta>\half$, thus beating the simple bound coming from the prime orbit theorem.
In a nutshell, we have reached the lower bound (for all $\varepsilon>0$), 
$$S(p)\geq C_\varepsilon (p-1)e^{(2\delta-1/2-\varepsilon)T}.$$
Keeping that lower bound in mind, we now turn to upper bounds using Proposition \ref{ExplicitF}.
Writing
$$S(p)=\vert I(id,T)\vert^2 +\sum_{\rr\neq id} \vert I(\rr,T)\vert^2, $$
and using the bound (\ref{trivial}) combined with the conclusion of Proposition \ref{ExplicitF}, we get 
$$S(p)=O(e^{(2\delta+\varepsilon)T})+O\left (  \sum_{\rr\neq id} d_\rr^2 (\log(d_\rr))^2 e^{2(\sigma+\varepsilon)T}  \right).$$
Using the formula
$$\vert \G \vert =\sum_{\rr} d_\rr^2,$$
combined with the fact that $\vert \G \vert=p(p^2-1)=O(p^3)$, 
we end up with
$$S(p)=O(e^{(2\delta-\varepsilon)T})+O\left (p^3\log(p) e^{2(\sigma+\varepsilon)T}  \right). $$
Since $T=\beta \log(p)$, we have obtained for all $p$ large \footnote{Note that the $\log(p)$ term has been absorbed in $p^\varepsilon$.}
$$Cp^{(2\delta-1/2-\varepsilon)\beta}\leq p^{(2\delta+\varepsilon)\beta-1}+p^{2+2(\sigma+\varepsilon)\beta+\varepsilon}.$$ 
Remark that since $\beta<2$, then if $\varepsilon$ is small enough we always have
$$(2\delta+\varepsilon)\beta-1< (2\delta-1/2-\varepsilon)\beta,$$
so up to a change of constant $C$, we actually have for all large $p$
$$Cp^{(2\delta-1/2-\varepsilon)\beta}\leq p^{2+2(\sigma+\varepsilon)\beta+\varepsilon}.$$
We have contradiction for $p$ large provided
$$\sigma< (\delta-\frac{1}{4}-\frac{1}{\beta})-\varepsilon-\frac{\varepsilon}{2\beta}.$$
Since $\beta$ can be taken arbitrarily close to $2$ and $\varepsilon$ arbitrarily close to $0$, we have a contradiction whenever
$\delta>\frac{3}{4}$ and $\sigma<\delta-\frac{3}{4}$. Therefore for all $p$ large, at least one of the $L$-function $L_\Gamma(s,\rr)$ 
for non trivial $\rr$ has to vanish inside the rectangle
$$\left \{ \delta-\textstyle{\frac{3}{4}}-\epsilon\leq \Re(s)\leq \delta\ \mathrm{and}\ 
\vert \Im(s) \vert \leq \left (\log(\log(p))\right)^{1+\alpha}\right \},$$
but then by the product formula we know that this zero appear as a zero of $Z_{\Gamma(p)}(s)$ with multiplicity $d_\rho$ which is greater or equal to $\frac{p-1}{2}$ by Frobenius. The main theorem is proved. $\square$

\bigskip
We end by a few comments. It is rather clear to us that this strategy should work without major modification for congruence subgroups
of arithmetic groups arising from quaternion algebras and also in higher dimension i. e. for convex co-compact subgroups of 
$SL_2(\Z(i))$.
What is less clear is the possibility to obtain similar results for more general families of Galois covers since we have used arithmeticity
in a rather fundamental way (via exponential multiplicities in the length spectrum). 

It would be interesting to know if the $\log^{1+\epsilon}(\log(p))$ bound can be improved to a uniform constant. However, it would likely
require a completely different approach since $\log(\log(p))$ is the very limit one can achieve with compactly supported test functions.
Indeed, to achieve a uniform bound with our approach would require the use of test functions $\varphi \not \equiv 0$ with Fourier
bounds 
$$\vert \widehat{\varphi}(\xi)\vert \leq C_1 e^{\vert \Im(\xi)\vert} e^{-C_2\vert \Re(\xi)\vert},$$ 
but an application of the Paley-Wiener theorem shows that these test functions do not exist (they would be both compactly supported
and analytic on the real line). 

 \end{document}